\documentstyle[12pt]{article}
\begin{document}
\parskip 4pt

\title{Free boson representation of $
DY_{\hbar }\left( \widehat{sl}
\left( M+1|N+1\right) \right) $ at level one}
\author{{Bo-yu Hou$^{a}$\thanks{E-mail address: 
byhou@phy.nwu.edu.cn}, 
Wen-Li Yang$^{b}$\thanks{E-mail address: 
wlyang@th.physik.uni-bonn.de}
and Yi Zhen$^{a}$\thanks{E-mail address: zheny@phy.nwu.edu.cn}
} \\
{\small $~^{a}$Institute of Modern Physics , Northwest University
 Xian 710069 , P.R.China}\\
{\small $~^{b}$ Physikalishes Institut der Universit\"at Bonn, 
Nussalle 12, 53115 Bonn, Germany}
}
\maketitle

\begin{abstract} We construct a realization of the central
extension of super-Yangian double $
DY_{\hbar }\left( \widehat{sl}\left( M+1|N+1\right) \right) $ at 
level-one in terms of free boson fields with a continuous parameter.

\vspace{1.2cm}

\end{abstract}


\def\a{\alpha}
\def\b{\beta}
\def\d{\delta}
\def\e{\epsilon}
\def\ve{\varepsilon}
\def\g{\gamma}
\def\k{\kappa}
\def\l{\lambda}
\def\h{\frac{1}{2\hbar}}
\def\o{\omega}
\def\t{\theta}
\def\s{\sigma}
\def\D{\Delta}
\def\L{\Lambda}

\def\R{\overline{R}}
\def\S{\overline{S}}
\def\Sl{{sl(n)}}
\def\DR{{R_{\xi}^{DY}}}
\def\VR{{R^{Y}}} 
\def\hS{{\widehat{sl(n)}}}
\def\hG{{\widehat{gl(N|N)}}}
\def\R{{\cal R}}
\def\hR{{\hat{\cal R}}}
\def\C{{\cal C}}
\def\P{{\bf P}}
\def\Z2{{{\bf Z}_2}}
\def\Z{{\bf Z}_n}
\def\T{{\cal T}}
\def\H{{\cal H}}
\def\F{{\cal F}}
\def\V{\overline{V}}
\def\trho{{\tilde{\rho}}}
\def\tphi{{\tilde{\phi}}}
\def\tT{{\tilde{\cal T}}}
\def\uqsnh{{U_q[\widehat{sl}(M+1|N+1)]}}
\def\uqgnh{{U_q[\widehat{gl}(N|N)]}}
\def\uq1h{{U_q[\widehat{gl(1|1)}]}}
\def\uqg2h{{U_q[\widehat{gl(2|2)}]}}
\def\ady{{DY_{\hbar }\left( \widehat{sl}
\left( M+1|N+1\right) \right)}}
\def\dy{{DY_{\hbar }\left( sl
\left( M+1|N+1\right) \right)}}


\def\beq{\begin{equation}}
\def\eeq{\end{equation}}
\def\bea{\begin{eqnarray}}
\def\eea{\end{eqnarray}}
\def\ba{\begin{array}}
\def\ea{\end{array}}
\def\no{\nonumber}
\def\lt{\left}
\def\rt{\right}
\newcommand{\bq}{\begin{quote}}
\newcommand{\eq}{\end{quote}}

\newtheorem{Theorem}{Theorem}
\newtheorem{Definition}{Definition}
\newtheorem{Proposition}{Proposition}
\newtheorem{Lemma}{Lemma}
\newtheorem{Corollary}{Corollary}
\newcommand{\proof}[1]{{\bf Proof. }
        #1\begin{flushright}$\Box$\end{flushright}}

\newcommand{\sect}[1]{\setcounter{equation}{0}\section{#1}}
\renewcommand{\theequation}{\thesection.\arabic{equation}}

\section{Introduction}

The algebraic analysis method based on infinite-dimensional highest 
weight representations of $non-abelian$ symmetries such as deformed 
infinite-dimensional (super) algebras has proved eminently 
successful in solving lattice integrable models\cite{Jim94,Koy94,Hou99,
Yan99} and completely integrable field theories\cite{Sim92,Kon97}. 
A powerful approach for studying the highest weight representations 
is the bosonization technique\cite{Fre88,Ber89} which allows one to 
explicitly construct these objects in terms of the deformed free 
bosonic fields. 

Free bosonic realization of level-one representations has been 
constructed for most quantum affine algebras\cite{Fre88,Ber89,Jin97}. This 
bosonization technique has also been extended to type I quantum affine 
super algebras $\uqsnh$, $M\neq N$\cite{Kim97}, $\uqgnh$\cite{Zha99} and 
$U_q[osp(2|2)^{(2)}]$\cite{Yan991}. On the other hand, free bosonic
realizations of  Yangian-deformed algebras have been constructed for 
$DY_{\hbar}(\widehat{sl}_2)$ at level $k$\cite{Kon971}, 
$DY_{\hbar}(\widehat{gl}_N)$ at level one\cite{Ioh96} and level
$k$\cite{Din98}, and 
$A_{\hbar,\eta}(\widehat{sl}_2)$ at level one\cite{Kho98}. It should be  
remarked that the free bosonic field with  continuous parameter initiated 
by Jimbo et al\cite{Jim96} is 
powerful in constructing the representations of Yangian-deformed 
algebras\cite{Hou97,Din98,Kho98}. The purpose of this paper is to
construct  the central extension super-Yangian double
(super-affine-Yangian algebra)   $\ady$  algebra at level one in terms of
free bosonic fields 
with continuous parameter.

This paper  is arranged as follows. In section 2 we give the defining
relations of  Drinfield currents of super-affine-Yangian algebra   
$\ady$. The free bosonic realization of
$\ady$ are constructed in section 3. Appendix contains some detailed
calculations. 

\sect{Drinfeld currents of $\ady$}

We will study the super-affine-Yangian algebra  $\ady$ for 
$M,N=0,1,\cdots$. The super-affine Yangian algebra $\ady$
is an associative algebra  generated by the  Drinfeld currents: $X_{i}^{\pm
}\left( \mu \right) $, $\Psi _{j}^{\pm }\left( \mu \right)$, $
i=1,2,...M+N+1,j=1,2,...,M+N$ and a central element $c$. The generators  
$\Psi _{j}^{\pm }\left(\mu \right) $ are invertible. The ${\bf
Z}_2$-grading of
the generators is : $\left[
X_{M+1}^{\pm }\left( \mu \right) \right] =1$ and zero otherwise . The
defining relations are \cite{Zha97,Cai97} 
\begin{eqnarray*}
\Psi_{i+1}^{-}( \nu ) ^{-1}\Psi_{i}^{+}( \mu) &=&\frac{(\mu-\nu-2\hbar)
(\mu-\nu +2\hbar) }{(\mu-\nu)(\mu-\nu)}\Psi_{i}^{+}(\mu)
\Psi _{i+1}^{-}(\nu) ^{-1} ,~~ M+N+1\geq i\geq M+1, \\
\Psi _{i+1}^{-}(\nu)^{-1}\Psi _{i}^{+}(\mu)
 &=&\frac{(\mu-\nu)(\mu -\nu) }
         {(\mu-\nu-2\hbar)(\mu-\nu+2\hbar) }
   \Psi_{i}^{+}(\mu)\Psi _{i+1}^{-}(\nu) ^{-1},\; \quad i<M+1,\\
\Psi _{i+1}^{+}(\nu) ^{-1}\Psi _{i}^{-}(\mu) &=&
\frac{(\mu -\nu -2\hbar) (\mu -\nu +2\hbar) }
{(\mu -\nu)(\mu -\nu)}\Psi _{i}^{-}( \mu
) \Psi _{i+1}^{+}( \nu) ^{-1},\;\quad i<M+1, \\
\Psi _{i+1}^{+}(\nu) ^{-1}\Psi _{i}^{-}(\mu) &=&
\frac{(\mu -\nu)(\mu -\nu)}{( \mu -\nu-2\hbar)(\mu-\nu +2\hbar)}\Psi _{i}^{-}
( \mu
) \Psi _{i+1}^{+}( \nu) ^{-1},\;\quad M+N+1\geq i\geq M+1, \\
\Psi _{i}^{+}(\mu) \Psi _{i}^{\_}(\nu) &=&\frac{
(\mu -\nu -3\hbar )( \mu -\nu +3\hbar ) }{( \mu
-\nu +\hbar )( \mu -\nu -\hbar ) }\Psi _{i}^{-}( \nu) \Psi _{i}^{+}( \mu
),\;\quad i<M+1, \\
\Psi _{i}^{+}\left( \mu \right) \Psi _{i}^{\_}\left( \nu \right) &=&\frac{
\left( \mu -\nu -\hbar \right) \left( \mu -\nu +\hbar \right) }{\left( \mu
-\nu +3\hbar \right) \left( \mu -\nu -3\hbar \right) }\Psi _{i}^{-}\left(
\nu \right) \Psi _{i}^{+}\left( \mu \right),\;\quad M+N+1\geq i>M+1 \\
\Psi _{M+1}^{+}\left( \mu \right) \Psi _{M+1}^{\_}\left( \nu \right) &=&\Psi
_{M+1}^{-}\left( \nu \right) \Psi _{M+1}^{+}\left( \mu \right) \\
\end{eqnarray*}
\vspace {-1.2cm}
\begin{eqnarray*}
\Psi _{i}^{+}\left( \mu \right) ^{-1}X_{i}^{+}\left( \nu \right) \Psi
_{i}^{+}\left( \mu \right) &=&\frac{\left( \mu _{+}-\nu -3\hbar \right) }{%
\left( \mu _{+}-\nu +\hbar \right) }\;X_{i}^{+}\left( \nu \right) 
,\;\quad M+N+1\geq i>M+1, \\
\Psi _{i}^{+}\left( \mu \right) ^{-1}X_{i}^{+}\left( \nu \right) \Psi
_{i}^{+}\left( \mu \right) &=&\frac{\left( \mu _{+}-\nu +\hbar \right) }{%
\left( \mu _{+}-\nu -3\hbar \right) }\;X_{i}^{+}\left( \nu \right) 
,\;\quad i<M+1, \\
\Psi _{i}^{-}\left( \mu \right) ^{-1}X_{i}^{+}\left( \nu \right) \Psi
_{i}^{-}\left( \mu \right) &=&\frac{\left( \mu _{-}-\nu +3\hbar \right) }{%
\left( \mu _{-}-\nu -\hbar \right) }\; X_{i}^{+}\left( \nu \right)
,\;\quad i<M+1, \\
\Psi _{i}^{-}\left( \mu \right) ^{-1}X_{i}^{+}\left( \nu \right) \Psi
_{i}^{-}\left( \mu \right) &=&\frac{\left( \mu _{-}-\nu -\hbar \right) }{%
\left( \mu _{-}-\nu +3\hbar \right) }\;X_{i}^{+}\left( \nu \right) 
,\;\ M+N+1\; \;\geq i>M+1, \\
\Psi _{i+1}^{+}\left( \mu \right) ^{-1}X_{i}^{+}\left( \nu \right) \Psi
_{i+1}^{+}\left( \mu \right) &=&\frac{\left( \mu _{+}-\nu -2\hbar \right) }{%
\left( \mu _{+}-\nu \right) }\; X_{i}^{+}\left( \nu \right)
,\;\quad i<M+1, \\
\Psi _{i+1}^{+}\left( \mu \right) ^{-1}X_{i}^{+}\left( \nu \right) \Psi
_{i+1}^{+}\left( \mu \right) &=&\frac{\left( \mu _{+}-\nu \right) }{\left(
\mu _{+}-\nu -2\hbar \right) }\;X_{i}^{+}\left( \nu \right)
,\;\quad  M+N+1\geq i\geq M+1, \\
\Psi _{i+1}^{-}\left( \mu \right) ^{-1}X_{i}^{+}\left( \nu \right) \Psi
_{i+1}^{-}\left( \mu \right) &=&\frac{\left( \mu _{-}-\nu \right) }{\left(
\mu _{-}-\nu +2\hbar \right) }\; X_{i}^{+}(\nu),\;\quad i<M+1, \\
\Psi _{i+1}^{-}\left( \mu \right) ^{-1}X_{i}^{+}\left( \nu \right) \Psi
_{i+1}^{-}\left( \mu \right) &=&\frac{\left( \mu _{-}-\nu +2\hbar \right) }{%
\left( \mu _{-}-\nu \right) } \; X_{i}^{+}\left( \nu \right),\; \quad  
M+N+1\geq i\geq M+1, \\
\Psi _{i}^{+}\left( \mu \right) ^{-1}X_{i}^{-}\left( \nu \right) \Psi
_{i}^{+}\left( \mu \right) &=&\frac{\left( \mu _{-}-\nu -\hbar \right) }{%
\left( \mu _{-}-\nu +3\hbar \right) }\; X_{i}^{-}\left( \nu \right) 
,\;\quad i<M+1, \\
\Psi _{i}^{+}\left( \mu \right) ^{-1}X_{i}^{-}\left( \nu \right) \Psi
_{i}^{+}\left( \mu \right) &=&\frac{\left( \mu _{-}-\nu +3\hbar \right) }{%
\left( \mu _{-}-\nu -\hbar \right) }\; X_{i}^{-}\left( \nu
\right),\;\quad M+N+1\geq i>M+1, \\
\Psi _{i}^{-}\left( \mu \right) ^{-1}X_{i}^{-}\left( \nu \right) \Psi
_{i}^{-}\left( \mu \right) &=&\frac{\left( \mu _{+}-\nu -3\hbar \right) }{%
\left( \mu _{+}-\nu +\hbar \right) }\; X_{i}^{-}\left( \nu \right) 
,\;\quad i<M+1, \\
\Psi _{i}^{-}\left( \mu \right) ^{-1}X_{i}^{-}\left( \nu \right) \Psi
_{i}^{-}\left( \mu \right) &=&\frac{\left( \mu _{+}-\nu +\hbar \right) }{%
\left( \mu _{+}-\nu -3\hbar \right) }\; X_{i}^{-}\left( \nu \right) 
,\;\quad M+N+1\geq i>M+1, \\
\Psi _{i+1}^{+}\left( \mu \right) ^{-1}X_{i}^{-}\left( \nu \right) \Psi
_{i+1}^{+}\left( \mu \right) &=&\frac{\left( \mu _{-}-\nu +2\hbar \right) }{%
\left( \mu _{-}-\nu \right) }\; X_{i}^{-}\left( \nu \right) \quad i<M+1,
\\
\Psi _{i+1}^{+}\left( \mu \right) ^{-1}X_{i}^{-}\left( \nu \right) \Psi
_{i+1}^{+}\left( \mu \right) &=&\frac{\left( \mu _{-}-\nu \right) }{\left(
\mu _{-}-\nu +2\hbar \right) }\; X_{i}^{-}\left( \nu \right),\; \quad 
M+N+1\geq i\geq M+1, \\
\Psi _{i+1}^{-}\left( \mu \right) ^{-1}X_{i}^{-}\left( \nu \right) \Psi
_{i+1}^{-}\left( \mu \right) &=&\frac{\left( \mu _{+}-\nu \right) }{\left(
\mu _{+}-\nu -2\hbar \right) }\; X_{i}^{-}\left( \nu \right) ,\;\quad
i<M+1,\\
\Psi _{i+1}^{-}\left( \mu \right) ^{-1}X_{i}^{-}\left( \nu \right) \Psi
_{i+1}^{-}\left( \mu \right) &=&\frac{\left( \mu _{+}-\nu -2\hbar \right) }{%
\left( \mu _{+}-\nu \right) }\; X_{i}^{-}\left( \nu \right) ,\;\quad 
M+N+1\geq i\geq M+1, \\
\Psi _{M+1}^{\pm }\left( \mu \right)
^{-1}X_{M+1}^{\pm }\left( \nu \right) \Psi _{M+1}^{\pm }\left( \mu \right)
&=&X_{M+1}^{\pm }\left( \nu \right) \\
\Psi _{M+1}^{\pm }\left( \mu \right) ^{-1}X_{M+1}^{\mp }\left( \nu \right)
\Psi _{M+1}^{\pm }\left( \mu \right) &=&X_{M+1}^{\mp }\left( \nu
\right)\\
\left( \mu -\nu \mp 2\hbar \right) X_{i}^{\mp }\left( \mu \right) X_{i}^{\mp
}\left( \nu \right) &=&\left( \mu -\nu \pm 2\hbar \right) X_{i}^{\mp }\left(
\nu \right) X_{i}^{\mp }\left( \mu \right) ,\;\quad i<M+1, \\
\left( \mu -\nu \pm 2\hbar \right) X_{i}^{\mp }\left( \mu \right) X_{i}^{\mp
}\left( \nu \right) &=&\left( \mu -\nu \mp 2\hbar \right) X_{i}^{\mp }\left(
\nu \right) X_{i}^{\mp }\left( \mu \right) ,\;\quad M+N+1\geq i>M+1, \\
\left( \mu -\nu \pm \hbar \right) X_{i}^{\mp }\left( \mu \right)
X_{i+1}^{\mp }\left( \nu \right) &=&\left( \mu -\nu \mp \hbar \right)
X_{i+1}^{\mp }\left( \nu \right) X_{i}^{\mp }\left( \mu \right),\;\quad 
i<M+1, \\
\left( \mu -\nu \mp \hbar \right) X_{i}^{\mp }\left( \mu \right)
X_{i+1}^{\mp }\left( \nu \right) &=&\left( \mu -\nu \pm \hbar \right)
X_{i+1}^{\mp }\left( \nu \right) X_{i}^{\mp }\left( \mu \right)
,\;\quad M+N+1\geq i\geq M+1, 
\end{eqnarray*}
\vspace {-1.2cm}
\begin{eqnarray*}
\left[ X_{i}^{+}\left( \mu \right),~X_{j}^{-}\left(
\nu \right) \right] &=&-2\hbar \delta _{ij}\left( \delta \left( \mu _{-}-\nu
_{+}\right) \Psi _{i}^{-}\left( \nu _{+}\right) -\delta \left( \mu _{+}-\nu
_{-}\right) \Psi _{i}^{+}\left( \mu _{+}\right) \right),\;\quad 
i,j\neq M+1, \\
\left\{ X_{M+1}^{+}\left( \mu \right) ~,~ 
X_{M+1}^{-}\left( \nu \right) \right\} &=&2\hbar \left( \delta \left( \mu
_{-}-\nu _{+}\right) \Psi _{i}^{-}\left( \nu _{+}\right) -\delta \left( \mu
_{+}-\nu _{-}\right) \Psi _{i}^{+}\left( \mu _{+}\right) \right), \\
\end{eqnarray*}
where \lbrack X,Y\rbrack $\equiv XY-YX$ stands for a commutator and \{X, Y\}$%
\equiv XY+YX$ for \ an anti-commutator.
As for the serre relations for the Drinfeld currents, we will refer the
reader to the ref.\cite{Zha97}.

We should remark that our definition of $\ady$ is equivalent to that of 
Zhang's\cite{Zha97} with defining
\bea
\Psi _{j}^{\pm }( \mu)= k_{j+1}^{\mp
}(\mu) k_{j}^{\mp }(\mu )^{-1},\no
\eea
and  spectra-shifting 
\bea
\Psi_{j}^{\pm}(\mu)&\rightarrow& \Psi_{j}^{\pm}(\mu +(i-1)\hbar)
\quad i<M+1,\no\\
\Psi_{j}^{\pm}(\mu)& \rightarrow&
\Psi_{j}^{\pm}(\mu+(M-i)\hbar)\quad M+N+1\geq i\geq M+1.\no\\
X_{j}^{\pm}(\mu)&\rightarrow& X_{j}^{\pm}(\mu +(i-1)\hbar)
\quad i<M+1,\no\\
X_{j}^{\pm}(\mu)& \rightarrow&
X_{j}^{\pm}(\mu+(M-i)\hbar)\quad M+N+1\geq i\geq M+1.\no
\eea

\sect{Free Bosonic realization}
Now, we study the free bosonic realization of $\ady$ which give
us a level one representation. Free bosonic realizations of the level-one 
representations of quantum affine super algberas have been obtained for 
$\uqsnh$, $M\neq N$\cite{Kim97} and $\uqgnh$\cite{Zha99}. Our
representation can be regarded as extending  the above works to the
super-affine-Yangian case.  

Let introduce  bosonic oscillators 
$\{\widehat{a_i}(t),\widehat{b}_j(t),\widehat{c}_j(t),|i=1,\cdots,M+1, 
j=1,\cdots,N+1\}$ 
with a continuous
parameters $t ( t\in R-\left\{ 0\right\}$, which satisfy  the
following commutation relations
\bea
\left[ \widehat{a_{i}}\left( t\right) \widehat{a_{j}}\left(
t^{\prime }\right) \right] &=&\frac{1}{\left( i\hbar \right) ^{2}t}%
sh^{2}\left( i\hbar t\right) \delta _{ij}\delta \left( t+t^{\prime
}\right),\\
\left[ \widehat{b_{i}}\left( t\right)\widehat{b_{j}}\left(
t^{\prime }\right) \right] &=&-\frac{1}{\left( i\hbar \right) ^{2}t}%
sh^{2}\left( i\hbar t\right) \delta _{ij}\delta \left( t+t^{\prime
}\right),\\
\left[ \widehat{c_{i}}\left( t\right) \widehat{c_{j}}\left(
t^{\prime }\right) \right] &=&\frac{1}{\left( i\hbar \right) ^{2}t}%
sh^{2}\left( i\hbar t\right) \delta _{ij}\delta \left( t+t^{\prime
}\right),
\eea
and the other commutators vanish. 
Define bosonic oscillators $\widehat{\l_i}(t)(1\leq i\leq M+N+1)$ 
\bea
\widehat{\l}_i(t)&=&\widehat{a_i}(t)e^{i\hbar\frac{|t|}{2}}-
\widehat{a}_{i+1}(t)e^{-i\hbar\frac{|t|}{2}},~~i=1,\cdots,M,\no\\
\widehat{\l}_{M+1}(t)&=&\widehat{a}_{M+1}(t)e^{i\hbar\frac{|t|}{2}}+
\widehat{b}_{1}(t)e^{i\hbar\frac{|t|}{2}},\no\\
\widehat{\l}_{M+1+j}(t)&=&-\widehat{b_{j}}(t)e^{-i\hbar\frac{|t|}{2}}+
\widehat{b}_{j+1}(t)e^{i\hbar\frac{|t|}{2}},~~j=1,\cdots,N.\no
\eea
It is easy to verify that the bosonic oscillators satisfy 
\bea
\left[ \widehat{\lambda _{i}}\left( t\right) \widehat{%
\lambda _{j}}\left( t^{\prime }\right) \right] &=&\frac{1}{\left( i\hbar
\right) ^{2}t}sh\left( ia_{ij}\hbar t\right) sh\left( i\hbar t\right)
\delta
\left( t+t^{\prime }\right),
\eea
where $a_{ij}$ is the Cartan matrix of the  super algebra 
$sl(M+1|N+1)$ 
\bea
(a_{ij})=\left(\begin{array}{ccccccccc}
2&-1&&&&&&&\\
-1&2&\ddots&&&&&&\\
&\ddots&\ddots&-1&&&&&\\
&&-1&2&-1&&&&\\
&&&-1&0&1&&&\\
&&&&1&-2&\ddots&&\\
&&&&&\ddots&\ddots&1&\\
&&&&&&1&-2&1\\
&&&&&&&1&-2
\end{array}
\right)
~(1\leq i,j\leq M+N+1).\label{Cart}
\eea
Moreover, we introduce the zero mode
operators 
$\{a_i(0),b_j(0),c_j(0),Q_{a_i},Q_{b_j},Q_{c_j}|i=1,\cdots,M+1,
j=1,\cdots,N+1\}$ 
satisfying the following commutation relations
\bea
\left[ a_{i}( 0), Q_{a_{i}}\right] &=&\delta _{ij}, \\
\left[b_{i}( 0), Q_{b_{i}}\right] &=&-\delta _{ij}, \\
\left[c_{i}( 0), Q_{c_{i}}\right] &=&\delta _{ij}.
\eea
Define $Q_{\widehat{\l}_i}=Q_{a_i}-Q_{a_{i+1}}$ for $i=1,\cdots,M$, 
$Q_{\widehat{\l}_{M+1}}=Q_{a_{M+1}}+Q_{b_{1}}$ and 
$Q_{\widehat{\l}_{M+1+j}}=-Q_{b_j}+Q_{b_{i+1}}$ for $j=1,\cdots,N$.
Let us introduce the notations
\bea
\chi _{i}( \mu ;\beta )&=&\exp \left\{ -i\hbar \int_{-\infty
}^{0}\frac{dt}{sh\left( i\hbar t\right) }\widehat{\chi _{i}}\left( t\right)
e^{-i\hbar \beta t}e^{-i\mu t}\right\}\no\\
&&~~\times\exp \left\{ -i\hbar
\int_{-\infty }^{0}\frac{dt}{sh\left( i\hbar t\right) }\widehat{\chi _{i}}
\left( t\right) e^{-i\hbar \beta t}e^{-i\mu t}\right\} e^{Q_{\widehat{\chi
_{i}}}}, \\
\chi _{i}^{+}\left( \mu \right) &=&\exp \left\{ -2\hbar \int_{0}^{\infty }%
\widehat{\chi _{i}}\left( t\right) e^{-i\mu t}dt\right\},\\
\chi _{i}^{-}\left( \mu \right) &=&\exp \left\{ 2\hbar \int_{-\infty }^{0}%
\widehat{\chi _{i}}\left( t\right) e^{-i\mu t}dt\right\}.
\eea
In the following, we will adopt these notations for the other bosonic
fields $\{\l_i(\mu;\b),c_j(\mu;\b)$ $|i=1,\cdots,M+N+1,~j=1,\cdots,N+1\}$, 
for example, the bosonic field $c_j(\mu;\b)$ should be defined in the same 
way. We introduce  the $\hbar$-difference operator defined
by 
\bea 
D_{\hbar}\left( f(\beta )\right) =f(\beta -i\hbar )-f(\beta +i\hbar).
\eea

When calculating the normal order of the bosonic fields, one  often 
encounter an integral 
\[
\int_{0}^{\infty }F(t)dt,
\]
which is divergent at t=0 . Here we adopt the following 
regularization\cite{Jim96,Kho98}. The integral 
should be understood as the contour integral 
\[
\int_{C}F(t)\frac{\ln (-t)}{2\pi i}dt, 
\]
where $C$ indicate the circle from +$\infty $ to 0 in the upper half plain
and
0 to $\infty $ in the lower half plain

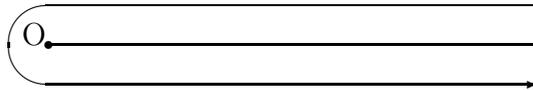
\begin{picture}(110,20)(10,20)
\put(165,10){\line(1,0){185}}
\put(165,10){\circle*{3}}
\put(155,10){O}
\put(350,10){\oval(400,30)[l]}
\put(340,-5){\vector(1,0){10}}
\end{picture} \\[6mm]

From the above regularization and after a straightforward calculation, 
we can obtain the normal order relations of the basic bosonic fields which 
will be given in the appendix A. Using the normal order relations, 
we can construct the level-one representation of $\ady$ in terms of free 
bosonic fields. Namely, 
the Drinfeld currents of $\ady$ at level-one are realized by the free
boson fields as follow
\bea
X_{i}^{+}\left( \mu \right) &=&:\lambda _{i}\left( \mu ;\frac{1}{2}\right)
D_{\hbar }\left[ c_{i-M-1}^{-1}\left( \mu ;0\right) \right] c_{i-M}\left(
\mu ;0\right) :,~~ M+N+1\geq i>M+1,\no\\ \\
X_{i}^{+}\left( \mu \right) &=&\frac{\sqrt{2\hbar }}{\exp \left( -\frac{1}{2}%
\gamma \right) }:\lambda _{i}\left( \mu ;\frac{1}{2}\right) c_{i-M}\left(
\mu ;0\right) : \prod_{j=1}^{M}\exp \left( -i\pi a_{j}\left(
0\right) \right),~~i=M+1,\no\\ \\
X_{i}^{+}( \mu ) &=&\frac{2\hbar }{\exp \left( -\gamma \right) }%
:\lambda _{i}\left( \mu ;\frac{1}{2}\right) :\exp \left( i\pi a_{i}
\left( 0\right) \right),~~i<M+1, \\
X_{i}^{-}\left( \mu \right) &=&:\lambda _{i}^{-1}\left( \mu ;-\frac{1}{2}%
\right) c_{i-M-1}^{-1}\left( \mu ;0\right) D_{\hbar }\lbrack
c_{i-M}^{-1}\left( \mu ;0\right) \rbrack :\quad M+N+1\geq i>M+1,\no\\ \\
X_{i}^{-}\left( \mu \right) &=&\frac{\sqrt{2\hbar }}{\exp \left( -\frac{1}{2}%
\gamma \right) }:\lambda _{i}^{-1}\left( \mu ;-\frac{1}{2}\right) D_{\hbar
}\lbrack c_{i-M}^{-1}\left( \mu ;0\right) \rbrack :\prod_{j=1}^{M}\exp
\left( i\pi a_{j}\left( 0\right) \right) \quad i=M+1,\no\\ \\
X_{i}^{-}\left( \mu \right) &=&\frac{2\hbar }{\exp \left( -\gamma \right) }%
:\lambda _{i}^{-1}\left( \mu ;-\frac{1}{2}\right) :\exp \left( -i\pi 
a_{i}\left( 0\right) \right) \quad i<M+1, \\
\Psi _{i}^{\pm }\left( \mu \right) &=&\lambda _{i}^{\pm }\left(
\mu\right).
\eea
This bosonic realization may help us to construct the bosonization of the 
supersymmetric $t-J$ model.
\vspace{2cm}
\section*{Acknownledgements}
W.-L. Yang would like to thank Prof. von Gehlen and the 
theoretical group of the Physikalishes Institut der 
Universit\"at Bonn for their kind hospitality. This 
work has been partly supported by the National Natural 
Science  Foundation of China. W.-L. Yang is supported by the 
Alexander von Humboldt Foundation. 
\vspace{2cm}
\section*{ Appendix A}

In this appendix, we give the normal order relations of the fundamental
bosonic fields. In order to calculate the normal order relations,
the following formula\cite{Kho98} is helpful 

\[
\int_{C}\frac{d\lambda \ln (-\lambda )}{2\pi i\lambda }\frac{e^{-x\lambda }}{%
1-e^{-\lambda /\eta }}=\ln \Gamma (\eta x)+(\eta x-\frac{1}{2})(\gamma -\ln
\eta )-\frac{1}{2}\ln 2\pi 
\]
where $\gamma$ is the Euler's constant.
Using the integral representation of  $\Gamma $
-functions and the definition of Drinfeld currents of $\ady$
at level-one, we can derive  
\begin{eqnarray*}
\l_{i}( \mu;\b_1 ) \l_{i}( \nu;\b_2) &=&
\frac{\Gamma \left[ 
\frac{1}{2\hbar }\left( \mu -\nu \right)
-\frac{(\b_1+\b_2)}{2}+
\frac{(1+a_{ii})}{2}\right] }{\Gamma
\left[ \frac{1}{2\hbar }\left( \mu -\nu \right) 
-\frac{(\b_1+\b_2)}{2}+
\frac{(1-a_{ii})}{2}
\right] }
~\frac{ e^{a_{ii}\gamma}}{\eta^{a_{ii}}}
 :\l_{i}( \mu;\b_1) \l_{i}( \nu;\b_2 ):, \\
\l_{i}( \mu;\b_1 ) \l_{i+1}( \nu;\b_2)  &=&\frac
{\Gamma \left[ \frac{1}{
2\hbar }\left( \mu -\nu \right) 
-\frac{(\b_1+\b_2)}{2}+
\frac{(1+a_{i,i+1})}{2}\right]}
{\Gamma
\left[ \frac{1}{2\hbar }\left( \mu -\nu \right) 
-\frac{(\b_1+\b_2)}{2}+
\frac{(1-a_{i,i+1})}{2}\right]}
~\frac{e^{a_{i,i+1}\g}}{\eta^{a_{i,i+1}}}:
\l_{i}\left( \mu;\b_1 \right)
\l_{i+1}\left( \nu;\b_2 \right) :, \\
\l_{i}( \mu;\b_1) \l_{j}( \nu;\b_2)
&=&
:\l_{i}( \mu;\b_1) \l_{j}( \nu;\b_2):,~~~|i-j|>1,\\       
c_{i}( \mu;0 ) c_{j}( \nu;0) &=&\d_{ij}~\frac{\Gamma
\left[ \frac{1}{2\hbar }\left( \mu -\nu \right) +1
\right] }{\Gamma \left[ 
\frac{1}{2\hbar }\left( \mu -\nu \right)
\right] }~
\frac{ e^{\gamma}}{\eta}
 :c_{i}( \mu;0) c_{j}( \nu;0 ):,
\end{eqnarray*}
where $\eta=\frac{1}{2\hbar}$ and $a_{ij}$ is the Cartan matrix
(\ref{Cart})

\end{document}